\documentclass[copyright,creativecommons]{eptcs}

\usepackage{setspace, amssymb, amsmath, amsfonts, amsthm, listings, epsfig,url}
\newtheorem{thm}{Theorem}
\newtheorem{defn}[thm]{Definition}
\newtheorem{lem}[thm]{Lemma}
\newtheorem{prop}[thm]{Proposition}
\newcommand{\comment}[1]{}
\newcommand{\DFT}[1]{\hat{#1}}
\newcommand{\IDFT}[1]{\check{#1}}
\newcommand{\mymod}[1]{\mathbin{~\mathrm{mod}}~#1}
\newcommand{\remove}[1]{}

\providecommand{\event}{CCA 2010} % Name of the event you are submitting to
\usepackage{breakurl}             % Not needed if you use pdflatex only.

\title{A Rigorous Extension of the Sch\"onhage-Strassen Integer Multiplication Algorithm Using Complex Interval Arithmetic}

\author{Thomas Steinke 
\institute{Department of Mathematics and Statistics\\
University of Canterbury\\%\thanks{A fine university.}\\
Christchurch, New Zealand}
\email{tas74@student.canterbury.ac.nz}
\and Raazesh Sainudiin
\institute{Department of Mathematics and Statistics\\
University of Canterbury\\%\thanks{A fine university.}\\
Christchurch, New Zealand}
\email{r.sainudiin@math.canterbury.ac.nz}}

\def\titlerunning{Integer Multiplication Using Interval Arithmetic}
\def\authorrunning{T.~Steinke \& R.~Sainudiin}
\begin{document}
\maketitle

\begin{abstract}
Multiplication of $n$-digit integers by long multiplication requires $O(n^2)$ operations and can be time-consuming. In 1970 A.~Sch\"onhage and V.~Strassen published an algorithm capable of performing the task with only $O(n \log(n))$ arithmetic operations over $\mathbb{C}$; naturally, finite-precision approximations to $\mathbb{C}$ are used and rounding errors need to be accounted for. Overall, using variable-precision fixed-point numbers, this results in an $O(n (\log(n))^{2 + \varepsilon})$-time algorithm. However, to make this algorithm more efficient and practical we need to make use of hardware-based floating-point numbers. How do we deal with rounding errors? and how do we determine the limits of the fixed-precision hardware? Our solution is to use interval arithmetic to guarantee the correctness of results and determine the hardware's limits. We examine the feasibility of this approach and are able to report that 75,000-digit base-256 integers can be handled using double-precision containment sets. This clearly demonstrates that our approach has practical potential; however, at this stage, our implementation does not yet compete with commercial ones, but we are able to demonstrate the feasibility of this technique.
\end{abstract}

\section{Introduction}

Multiplication of very large integers is a crucial subroutine of many algorithms such as the RSA cryptosystem \cite{RSA}.  Consequently, much effort has gone into finding fast and reliable multiplication algorithms; \cite{KNUTH} discusses several methods. The asymptotically-fastest known algorithm \cite{FUERER} requires $n \log(n) 2 ^ {O(\log^*(n))}$ steps, where $\log^*$ is the iterated logarithm --- defined as the number of times one has to repeatedly take the logarithm before the number is less than $1$. However, being asymptotically-faster does not translate to being faster in practice. We shall concern ourselves with the practicalities of the subject; we will analyse our algorithm's performance on a finite range of numbers.

The algorithm we are studying here is based on the first of two asymptotically fast multiplication algorithms by A.~Sch\"onhage and V.~Strassen \cite{SCHOENHAGE}. These algorithms are based on the convolution theorem and the fast Fourier transform. The first algorithm (the one we are studying) performs the discrete Fourier transform over $\mathbb{C}$ using finite-precision approximations. The second algorithm uses the same ideas as the first, but it works over the finite ring $\mathbb{Z}_{2 ^ {2 ^ n} + 1}$ rather than the uncountable field $\mathbb{C}$.  We wish to point out that ``\emph{the} Sch\"onhage-Strassen algorithm'' usually refers to the second algorithm. However, in this document we use it to refer to the first $\mathbb{C}$-based algorithm.

From the theoretical viewpoint, the second algorithm is much nicer than the first. The second algorithm does not require the use of finite-precision approximations to $\mathbb{C}$.  Also, the second algorithm requires $O(n \log(n) \log(\log(n)))$ steps to multiply two $n$-bit numbers, making it asymptotically-faster than the first algorithm. However, the second algorithm is much more complicated than the first, and it is outperformed by asymptotically-slower algorithms, such as long multiplication, for small-to-medium input sizes. In practice, both of the Sch\"onhage-Strassen algorithms are rarely used.

The first Sch\"onhage-Strassen Algorithm is more elegant, if the finite-precision approximations are ignored. More importantly, it is faster in practice. Previous studies \cite{CFASTER} have shown that the first algorithm can be faster than even highly-optimised implementations of the second. However, the first algorithm's reliance on finite-precision approximations, despite exact answers being required, leads to it being discounted.

The saving grace of the Sch\"onhage-Strassen algorithm is that at the end of the computation an integral result will be obtained. So the finite-precision approximations are rounded to integers. Thus, as long as rounding errors are sufficiently small for the rounding to be correct, an exact answer will be obtained. Sh\"onhage and Strassen showed that fixed-point numbers with a variable precision of $O(\log(n))$ bits would be sufficient to achieve this.

For the Sch\"onhage-Strassen algorithm to be practical, we need to make use of hardware-based floating-point numbers; software-based variable-precision numbers are simply too slow. However, we need to be able to manage the rounding errors. At the very least, we must be able to detect when the error is too large and more precision is needed. The usual approach to this is to prove some kind of worst-case error bound (for an example, see \cite{ABOUND}). Then we can be sure that, for sufficiently small inputs, the algorithm will give correct results. However, worst-case bounds are rarely tight. We propose the use of dynamic error bounds using existing techniques from computer-aided proofs.

Dynamic error detection allows us to move beyond worst-case bounds. For example, using standard single-precision floating-point numbers, our na\"ive implementation of the Sch\"onhage-Strassen algorithm sometimes gave an incorrect result when we tried multiplying two 120-digit base-256 numbers, but it usually gave correct results. Note that by a `na\"ive implementation' we simply mean a modificiation of the Sch\"onhage-Strassen algorithm that uses fixed-precison floating-point arithmetic and does not guarantee correctness. A worst-case bound would not allow us to use the algorithm in this case, despite it usually being correct. Dynamic error detection, however, would allow us to try the algorithm, and, in the rare instances where errors occur, it would inform us that we need to use more precision.

We will use complex interval containment sets for all complex arithmetic operations. This means that at the end of the computation, where we would ordinarily round to the nearest integer, we simply choose the unique integer in the containment set. If the containment set contains multiple integers, then we report an error. This rigorous extension of the Sch\"onhage-Strassen algorithm therefore constitutes a computer-aided proof of the desired product. When an error is detected, we must increase the precision being used or we must use a different algorithm.

For those unfamiliar with the Sch\"onhage-Strassen algorithm or with interval arithmetic, we describe these in section \ref{SectionAlg}. Then, in section \ref{SectionResults}, we show the empirical results of our study. Section \ref{Conclusion}, our conclusion, briefly discusses the implications of our results.

\section{The Algorithm} \label{SectionAlg}

For the sake of completeness we explain the Sch\"onhage-Strassen algorithm, as it is presented in \cite{SCHOENHAGE}. We also explain how we have modified the algorithm using interval arithmetic in subsection \ref{IntervalArithmetic}. Those already familiar with the material may skip all or part of this section.

We make the convention that a positive integer $x$ is represented in base $b$ (usually $b = 2^k$ for some $k \in \mathbb{N}$) as a vector ${x} \in \mathbb{Z}_b^n := \{0,1,2,\ldots,b-1\}^n$; the value of $x$ is
$$x = \sum_{i = 0}^{n - 1} x_i b ^ i.$$

\subsection{Basic Multiplication Algorithm}
The above definition immediately leads to a formula for multiplication. Let $x$ and $y$ be positive integers with representations ${x} \in \mathbb{Z}_b ^ n$ and ${y} \in \mathbb{Z}_b ^ m$. Then
$$x y = \left( \sum_{i = 0}^{n - 1} x_i b ^ i \right) \left( \sum_{j = 0}^{m - 1} y_j b ^ j \right) = \sum_{i = 0}^{n + m - 2} \sum_{j = \max\{ 0, i - m + 1 \}}^{\min\{ n - 1, i \}} x_j y_{i - j} b ^ {i} = \sum_{i = 0}^{n + m - 1} z_i b ^ {i}.$$
Of course, we cannot simply set $z_i = \sum_{j = \max\{ 0, i - m + 1 \}}^{\min\{ n - 1, i \}} x_j y_{i - j}$; this may violate the constraint that $0 \leq {z}_i \leq b - 1$ for every $i$. We must `carry' the `overflow'. This leads to the long multiplication algorithm (see \cite{KNUTH}).
\\\parbox{200pt}{
\begin{singlespace}
\texttt{
\begin{tabbing}
\hspace{1cm}\=$~~~~~$\=$~~~~~$\=$~~~~~$\=$~~~~~$\=$~~~~~$\=\\
\>{\bf The Long Multiplication Algorithm}\\
\>1.\>Input: ${x} \in \mathbb{Z}_b ^ n$ and ${y} \in \mathbb{Z}_b ^ m$\\
\>2.\>Output: ${z} \in \mathbb{Z}_b ^ {n + m}$ \# $z = x y$\\
\>3.\>Set $c = 0$. \# $c$ = carry\\
\>4.\>For $i = 0$ up to $n + m - 1$ do \{\\
\>5.\>\>Set $s = 0$. \# $s$ = sum\\
\>6.\>\>For $j = \max\{ 0, i - m + 1 \}$ up to $\min\{ n - 1, i \}$ do \{\\
\>7.\>\>\>Set $s = s + {x}_j {y}_{i - j}$.\\
\>8.\>\>\}.\\
\>9.\>\>Set ${z}_i = (s + c) \mymod{b}$.\\
\>10.\>\>Set $c = \lfloor (s + c) / b \rfloor$.\\
\>11.\>\}.\\
\>12.\>\# $c = 0$ at the end.
\end{tabbing}
}
\end{singlespace}
}\\
This algorithm requires $O(m n)$ steps (for a fixed $b$).  Close inspection of the long multiplication algorithm might suggest that $O(m n \log(\min\{m, n\}))$ steps are required as the sum $s$ can become very large. However, adding a bounded number (${x}_j {y}_{i - 1} < b^2$) to an unbounded number ($s$) is, \emph{on average}, a constant-time operation.

\subsection{The Discrete Fourier Transform} \label{SectionDFT}

The basis of the Sch\"onhage-Strassen algorithm is the discrete Fourier transform and the convolution theorem. The discrete Fourier transform is a map from $\mathbb{C}^n$ to $\mathbb{C}^n$. In this section we will define the discrete Fourier transform and we will show how it and its inverse can be calculated with $O(n \log(n))$ complex additions and multiplications.  See \cite{KNUTH} for further details.

\begin{defn}[Discrete Fourier Transform]
Let $x \in \mathbb{C}^n$ and let $\omega := e ^{\frac{2 \pi {\rm i}}{n}}$. Then define the discrete Fourier transform $\DFT{x} \in \mathbb{C}^n$ of $x$ by
$$\DFT{x}_i := \sum_{j = 0}^{n - 1} x_j \omega ^ {i j} ~~ (0 \leq i \leq n - 1).$$
\end{defn}

There is nothing special about our choice of $\omega$; the consequences of the following lemma are all that we need $\omega$ to satisfy. Any other element of $\mathbb{C}$ with the same properties would suffice.

\begin{lem} \label{LemmaOmegaProperties}
Let $n > 1$ and $\omega = e ^ {\frac{2 \pi {\rm i}}{n}}$. Then
$$\omega ^ n = 1 \text{ and } \omega ^ k \ne 1 \text{ for all } 0 < k < n$$
and, for all $0 < k < n$,
$$\sum_{i = 0}^{n - 1} \omega ^ {i k} = 0.$$
\end{lem}

Note that the case where $n = 1$ is uninteresting, as $\omega = 1$ and the discrete Fourier transform is the identity mapping $\DFT{x} = x$.

\remove{
\begin{proof}
Firstly, $$\omega ^ n = \left( e ^ {\frac{2 \pi {\rm i}}{n}} \right) ^ n = e ^ {2 \pi {\rm i}} = 1.$$ We know that $e ^ {\theta} = 1$ if and only if $\theta = 2 \pi {\rm i} m$ for some $m \in \mathbb{Z}$. Thus, if $\omega ^ k = 1$, then $k$ must be a multiple of $n$, which eliminates the possibility that $0 < k < n$.

Fix $k$ with $0 < k < n$ and let $s_k := \sum_{i = 0}^{n - 1} \omega ^ {i k}$. Then
$$\omega ^ k s_k = \sum_{i = 0}^{n - 1} \omega ^ {(i + 1) k} = \sum_{i = 1}^{n} \omega ^ {i k} = \sum_{i = 1}^{n - 1} \omega ^ {i k} + \omega ^ {k n} = \sum_{i = 1}^{n - 1} \omega ^ {i k} + 1 = \sum_{i = 1}^{n - 1} \omega ^ {i k} + \omega ^ {0 k} = \sum_{i = 0}^{n - 1} \omega ^ {i k} = s_k.$$
So $\omega ^ k s_k = s_k$. If $s_k \ne 0$, then we can divide by $s_k$ to get $\omega ^ k = 1$, which is impossible. So $s_k = 0$.
\end{proof}
}

Now we can prove that the discrete Fourier transform is a bijection.

\begin{prop}[Inverse Discrete Fourier Transform]
Let $x \in \mathbb{C}^n$ and let $\omega = e ^{\frac{2 \pi {\rm i}}{n}}$. Define $\IDFT{x} \in \mathbb{C}^n$ by
$$\IDFT{x}_i := \frac{1}{n} \sum_{j = 0}^{n - 1} x_j \omega ^ {- i j} ~~ (0 \leq i \leq n - 1).$$
Then this defines the inverse of the discrete Fourier transform --- that is, if $y = \DFT{x}$, then $\IDFT{y} = x$.
\end{prop}

\remove{
\begin{proof}
Fix $x \in \mathbb{C}^n$, let $y = \DFT{x}$ and let $z = \IDFT{y}$. We wish to show that $z = x$. If $n = 1$, then this is trivial, as $x = y = z$, so we may assume that $n > 1$. First of all, it follows from Lemma \ref{LemmaOmegaProperties} that, if $l \in \mathbb{Z}$ and $n$ does not divide $l$, then
$$\sum_{i = 0}^{n - 1} \omega ^ {i l} = 0.$$
If, on the other hand, $n$ divides $l$, then
$$\sum_{i = 0}^{n - 1} \omega ^ {i l} = n.$$
Now, fixing $i$ with $0 \leq i \leq n - 1$, we have
\begin{eqnarray*}
z_i &=& \frac{1}{n} \sum_{j = 0}^{n - 1} y_j \omega ^ {- i j}
= \frac{1}{n} \sum_{j = 0}^{n - 1} \left( \sum_{k = 0}^{n - 1} x_k \omega ^ {j k} \right) \omega ^ {- i j}
= \frac{1}{n} \sum_{k = 0}^{n - 1} \sum_{j = 0}^{n - 1} x_k \omega ^ {j k} \omega ^ {- i j}
= \frac{1}{n} \sum_{k = 0}^{n - 1} x_k \sum_{j = 0}^{n - 1} \omega ^ {j (k - i)}\\
&=& \frac{1}{n} \sum_{k = 0}^{n - 1} x_k \left\{ \begin{array}{cl} n, & \text{if $n$ divides $k - i$} \\ 0, & \text{otherwise} \\ \end{array} \right\}
= \sum_{k = 0}^{n - 1} x_k \left\{ \begin{array}{cl} 1, & \text{if $k - i = 0$} \\ 0, & \text{otherwise} \\ \end{array} \right\}\\
&=& x_i.\\
\end{eqnarray*}
\end{proof}
}

Now we explain the fast Fourier transform; this is simply a fast algorithm for computing the discrete Fourier transform and its inverse.

Let $n$ be a power of $2$ and $x \in \mathbb{C}^n$ be given. Now define $x_\text{even}, x_\text{odd} \in \mathbb{C}^{n/2}$ by
\[
\left( x_\text{even} \right)_i = x_{2 i},\qquad
\left( x_\text{odd} \right)_i = x_{2 i + 1},
\]
for all $i$ with $0 \leq i \leq n / 2 - 1$.

Now the critical observation of the Cooley-Tukey fast Fourier transform algorithm is the following. Fix $i$ with $0 \leq i \leq n - 1$ and let $\omega = e ^ {\frac{2 \pi i}{n}}$. Then we have
\begin{eqnarray*}
\DFT{x}_i &=& \sum_{j = 0}^{n - 1} x_j \omega ^ {i j}\\
\remove{&=& \sum_{j = 0}^{n/2 - 1} x_{2 j} \omega ^ {2 i j} + \sum_{j = 0}^{n/2 - 1} x_{2 j + 1} \omega ^ {2 i j + i}\\}
&=& \sum_{j = 0}^{n/2 - 1} \left( x_\text{even} \right)_{j} \left( \omega ^ 2 \right) ^ {i j} + \omega ^ i \sum_{j = 0}^{n/2 - 1} \left( x_\text{odd} \right)_{j} \left( \omega ^ 2 \right) ^ {i j}\\
\remove{&=& \sum_{j = 0}^{n/2 - 1} \left( x_\text{even} \right)_{j} \left( \omega ^ 2 \right) ^ {(i \mymod{n / 2}) j} + \omega ^ i \sum_{j = 0}^{n/2 - 1} \left( x_\text{odd} \right)_{j} \left( \omega ^ 2 \right) ^ {(i \mymod{n / 2}) j}\\}
&=& \left( \DFT{x}_\text{even} \right)_{i \mymod{n / 2}} + \omega ^ i \left( \DFT{x}_\text{odd} \right)_{i \mymod{n / 2}}.\\
\end{eqnarray*}
Note that $\left( \omega^2 \right)^{n / 2} = 1$, so taking the modulus is justified. This observation leads to the following divide-and-conquer algorithm.
\\\parbox{200pt}{
\begin{singlespace}
\texttt{
\begin{tabbing}
\hspace{1cm}\=$~~~~~$\=$~~~~~$\=$~~~~~$\=$~~~~~$\=$~~~~~$\=\\
\>{\bf The Cooley-Tukey Fast Fourier Transform}\\
\>1.\>Input: $n = 2 ^ k$ and $x \in \mathbb{C} ^ n$\\
\>2.\>Output: $\DFT{x} \in \mathbb{C} ^ n$\\
\>3.\>function FFT($k$, $x$) \{\\
\>4.\>\>If $k = 0$, then $\DFT{x} = x$.\\
\>5.\>\>Partition $x$ into $x_\text{even}, x_\text{odd} \in \mathbb{C}^{n / 2}$.\\
\>6.\>\>Compute $\DFT{x}_\text{even}$ = FFT($k - 1$, $x_\text{even}$) by recursion.\\
\>7.\>\>Compute $\DFT{x}_\text{odd}$ = FFT($k - 1$, $x_\text{odd}$) by recursion.\\
\>8.\>\>Compute $\omega = e ^ {\frac{2 \pi {\rm i}}{n}}$.\\
\>9.\>\>For $i = 0$ up to $n - 1$ do \{\\
\>10.\>\>\>Set $\DFT{x}_i = \left( \DFT{x}_\text{even} \right)_{i \mymod{n / 2}} + \omega ^ i \left( \DFT{x}_\text{odd} \right)_{i \mymod{n / 2}}$.\\
\>11.\>\>\}.\\
\>12.\>\}.\\
\end{tabbing}
}
\end{singlespace}
}\\
It is easy to show that this algorithm requires $O(n \log(n))$ complex additions and multiplications. With very little modification we are also able to obtain a fast algorithm for computing the inverse discrete Fourier transform.

Note that, to compute $\omega$, we can use the recurrence
$$\omega_1 = 1, ~~ \omega_2 = -1, ~~ \omega_4 = {\rm i}, ~~ \omega_{2n} = \frac{1 + \omega_n}{\left| 1 + \omega_n \right|} ~~ (n \geq 3),$$
where $\omega_n = e ^ {\frac{2 \pi {\rm i}}{n}}$. Other efficient methods of computing $\omega$ are also available.

\subsection{The Convolution Theorem} \label{SectionConvolution}

We start by defining the convolution. Let $a, b \in \mathbb{C}^n$. We can interpret $a$ and $b$ as the coefficients of two polynomials --- that is,
$$f_a(z) = a_0 + a_1 z + \cdots + a_{n - 1} z ^ {n - 1}.$$
The convolution of $a$ and $b$ --- denoted by $a * b$ --- is, for our purposes, the vector of coefficients obtained by multiplying the polynomials $f_a$ and $f_b$. Thus we have $f_{a * b}(z) = f_a(z) f_b(z)$ for all $z \in \mathbb{C}$. Note that we can add `padding zeroes' to the end of the coefficient vectors without changing the corresponding polynomial.

The convolution theorem relates convolutions to Fourier transforms. We only use a restricted form.

\begin{thm}[Convolution Theorem]
Let $a, b \in \mathbb{C}^n$ and $c := a * b \in \mathbb{C}^m$, where $m = 2n - 1$. Pad $a$ and $b$ by setting
$$a^\prime = (a_0, a_1, \cdots, a_{n - 1}, 0, \cdots, 0), b^\prime = (b_0, b_1, \cdots, b_{n - 1}, 0, \cdots, 0) \in \mathbb{C}^m.$$
Then, for every $i$ with $0 \leq i \leq m - 1$,
$$\DFT{c}_i = \DFT{a^\prime}_i \DFT{b^\prime}_i.$$
\end{thm}

\remove{
\begin{proof}
By definition, if $\omega = e ^ {\frac{2 \pi i}{m}}$ and $0 \leq i \leq m - 1$, then
$$\DFT{c}_i = f_c(\omega ^ i) = f_a(\omega^i) f_b(\omega^i) = f_{a^\prime}(\omega^i) f_{b^\prime}(\omega^i) = \DFT{a^\prime}_i \DFT{b^\prime}_i.$$
\end{proof}
}

The convolution theorem gives us a fast method of computing convolutions and, thus, of multiplying polynomials. Given $a, b \in \mathbb{C}^n$, we can calculate $c = a * b$ using only $O(n \log(n))$ arithmetic operations as follows.
\comment{
\begin{itemize}
\item Let $k = \lceil \log_2(2 n - 1) \rceil$. (We need a sufficiently large power of two for the fast Fourier transform algorithm to work.)
\item First we pad $a$ and $b$ to get $a^\prime$ and $b^\prime$ in $\mathbb{C}^{2^k}$.
\item We calculate $\DFT{a^\prime}$ and $\DFT{b^\prime}$ using the fast Fourier transform.
\item We calculate $\DFT{c}$ using the convolution theorem --- that is, $\DFT{c}_i = \DFT{a^\prime}_i \DFT{b^\prime}_i$ ($0 \leq i \leq 2n - 2$).
\item We calculate $c$ from $\DFT{c}$ using the inverse fast Fourier transform.
\end{itemize}}
\\\parbox{200pt}{
\begin{singlespace}
\texttt{
\begin{tabbing}
\hspace{1cm}\=$~~~~~$\=$~~~~~$\=$~~~~~$\=$~~~~~$\=$~~~~~$\=\\
\>{\bf The Fast Convolution Algorithm}\\
\>1.\>Input: $a, b \in \mathbb{C}^n$\\
\>2.\>Output: $c = a * b \in \mathbb{C}^{m}$\\
\>3.\>Set $k = \lceil \log_2(2 n - 1) \rceil$ and $m = 2 ^ k$.\\
\>4.\>\# Pad $a$ and $b$ so they are in $\mathbb{C}^n$.\\
\>5.\>Set $a^\prime = (a_0, a_1, \cdots, a_{n - 1}, 0, \cdots, 0), b^\prime = (b_0, b_1, \cdots, b_{n - 1}, 0, \cdots, 0) \in \mathbb{C}^m$.\\
\>6.\>Compute $\DFT{a^\prime} = \text{FFT}(k, a^\prime)$ and $\DFT{b^\prime} = \text{FFT}(k, b^\prime)$.\\
\>7.\>For $0 \leq i \leq m - 1$, set $\DFT{c}_i = \DFT{a^\prime}_i \DFT{b^\prime}_i$.\\
\>8.\>Compute $c = \text{FFT}^{-1}(k, \DFT{c})$.\\
\end{tabbing}
}
\end{singlespace}
}\\

\subsection{The Sch\"onhage-Strassen Algorithm} \label{SectionSSA}

The Sch\"onhage-Strassen algorithm multiplies two integers by convolving them and then performing carrys. Let two base-$b$ integer representations be ${x}$ and ${y}$. We consider the digits as the coefficients of two polynomials. Then $x = f_{{x}}(b)$, $y = f_{{y}}(b)$ and $$x y = f_{{x}}(b) f_{{y}}(b) = f_{{x} * {y}}(b).$$
So, to compute $x y$, we can first compute ${x} * {y}$ in $O(n \log(n))$ steps and then we can evaluate $f_{{x} * {y}}(b)$. The evaluation of $f_{{x} * {y}}(b)$ to yield an integer representation ${z}$ is simply the process of performing carrys.
\\\parbox{200pt}{
\begin{singlespace}
\texttt{
\begin{tabbing}
\hspace{1cm}\=$~~~~~$\=$~~~~~$\=$~~~~~$\=$~~~~~$\=$~~~~~$\=\\
\>{\bf The Sch\"onhage-Strassen Algorithm}\\
\>1.\>Input: ${x} \in \mathbb{Z}_b ^ n$ and ${y} \in \mathbb{Z}_b ^ n$\\
\>2.\>Output: ${z} \in \mathbb{Z}_b ^ {2n}$ \# $z = x y$\\
\>3.\>Compute ${x} * {y}$ using the fast convolution Algorithm.\\
\>4.\>Set $c = 0$. \#carry\\
\>5.\>For $i = 0$ up to $2n - 2$ do \{\\
\>6.\>\>Set ${z}_i = \left( \left({x} * {y}\right)_i + c \right) \mymod{b}$.\\
\>7.\>\>Set $c = \lfloor \left( \left({x} * {y}\right)_i + c \right) / b \rfloor$.\\
\>8.\>\}.\\
\>9.\>Set ${z}_{2n - 1} = c$.\\
\end{tabbing}
}
\end{singlespace}
}\\
Clearly the Sch\"onhage-Strassen algorithm performs the multiplication using $O(n \log(n))$ complex arithmetic operations.

When finite-precision complex arithmetic is done, rounding errors are introduced. However, this can be countered: We know that ${x} * {y}$ must be a vector of integers. As long as the rounding errors introduced are sufficiently small, we can round to the nearest integer and obtain the correct result. Sch\"onhage and Strassen \cite{SCHOENHAGE} proved that $O(\log(b n))$-bit floating point numbers give sufficient precision.

\subsection{Interval Arithmetic} \label{IntervalArithmetic}

Our rigorous extension of the algorithm uses containment sets.  %Broadly speaking, containment sets are machine-representable sets which are guaranteed to contain the correct value.  Containment sets allow us to bound numerical errors caused by finite precision.  Arithmetic with containment sets is well-defined.
By replacing all complex numbers with complex containment sets, we can modify the Sch\"onhage-Strassen algorithm to find a containment set of ${x} * {y}$; if the containment set only contains one integer-valued vector, then we can be certain that this is the correct value.  We have used rectangular containment sets of machine-representable floating-point intervals with directed rounding to guarantee the desired integer product.  A brief overview of the needed interval analysis \cite{Moore2009} is given next.

Let $\underline{x},\overline{x}$ be real numbers with $\underline{x} \leq \overline{x}$. Let $[\underline{x},\overline{x}] = \{x \in \mathbb{R}: \underline{x} \leq x \leq \overline{x}\}$ be a closed and bounded real interval and let the set of all such intervals be $\mathbb{IR} = \{[\underline{x},\overline{x}] : \underline{x} \leq \overline{x} \, ; \, \underline{x},\overline{x} \in \mathbb{R}\}$.  Note that $\mathbb{R} \subset \mathbb{IR}$ since we allow thin or punctual intervals with $\underline{x}=\overline{x}$.  If $\star$ is one of the arithmetic operators $+$, $-$, $\cdot$, $/$, we define arithmetic over operands in $\mathbb{IR}$ by $[\underline{a},\overline{a}] \star [\underline{b},\overline{b}] := \{a \star b: a \in [\underline{a},\overline{a}], b \in [\underline{b},\overline{b}]\}$, with the exception that $[\underline{a},\overline{a}] / [\underline{b},\overline{b}]$ is undefined if $0 \in [\underline{b},\overline{b}]$.  Due to continuity and monotonicity of the operations and compactness of the operands, arithmetic over $\mathbb{IR}$ is given by real arithmetic operations with the bounds:
\begin{alignat*}{2}
[\underline{a},\overline{a}]+[\underline{b},\overline{b}]&= [\underline{a}+\underline{b}, \overline{a}+\overline{b}] \\
[\underline{a},\overline{a}]-[\underline{b},\overline{b}]&=  [\underline{a}-\overline{b}, \overline{a}-\underline{b}] \\
[\underline{a},\overline{a}]\cdot[\underline{b},\overline{b}]&= [\min\{ \underline{a}\underline{b},  \underline{a}\overline{b}, \overline{a}\underline{b}, \overline{a}\overline{b} \}, \max\{ \underline{a}\underline{b},  \underline{a}\overline{b}, \overline{a}\underline{b}, \overline{a}\overline{b} \}] \\
[\underline{a},\overline{a}]/ [\underline{b},\overline{b}]&= [\underline{a},\overline{a}] \cdot [1/\overline{b}, 1/\underline{b}], \, \text{if } \, 0 \notin [\underline{b},\overline{b}] \enspace .
\end{alignat*}
In addition to the above elementary operations over elements in $\mathbb{IR}$, our algorithm requires us to contain the range of the square root function over elements in $\mathbb{IR} \cap [0,\infty)$.  Once again, due to the monotonicity of the square root function over non-negative reals it suffices to work with the real image of the bounds
$
\sqrt{[\underline{x},\overline{x}]} = [\sqrt{\underline{x}}, \sqrt{\overline{x}}], \, \text{if } \, 0 \leq \underline{x} % \enspace .
$.  To complete the requirements for our rigorous extension of the Sch\"onhage-Strassen algorithm we need to extend addition, multiplication and division by a non-zero integer to elements in
$$
\mathbb{IC} := \left\{ [\underline{z},\overline{z}]:=[\underline{z_1},\overline{z_1}] + {\rm i} [\underline{z_2},\overline{z_2}] : \, [\underline{z_1},\overline{z_1}], [\underline{z_2},\overline{z_2}] \in \mathbb{IR} \right\} .%\enspace .
$$
Interval arithmetic over $\mathbb{IR}$ naturally extends to $\mathbb{IC}$, the set of rectangular complex intervals.  Addition and subtraction over $[\underline{z},\overline{z}],[\underline{w},\overline{w}] \in \mathbb{IC}$ given by
$$[\underline{z},\overline{z}] \pm [\underline{w},\overline{w}]= ( [\underline{z_1},\overline{z_1}] \pm [\underline{w_1},\overline{w_1}] ) + {\rm i} ([\underline{z_2},\overline{z_2}] \pm [\underline{w_2},\overline{w_2}])$$ are sharp but not multiplication or division due to rectangular wrapping effects.  Complex interval multiplication and division of a complex interval by a non-negative integer can be contained with real interval multiplications given by
$$[\underline{z},\overline{z}] \cdot [\underline{w},\overline{w}]  = ([\underline{z_1},\overline{z_1}] \cdot [\underline{w_1},\overline{w_1}] - [\underline{z_2},\overline{z_2}] \cdot [\underline{w_2},\overline{w_2}]) + {\rm i} ([\underline{z_1},\overline{z_1}] \cdot [\underline{w_2},\overline{w_2}] + [\underline{z_2},\overline{z_2}] \cdot [\underline{w_1},\overline{w_1}]).$$
% \\
%[\underline{z},\overline{z}]\cdot[\underline{w},\overline{w}]&=  \\
%[\underline{a},\overline{a}]/ [\underline{b},\overline{b}]&= , \, \text{if } \, 0 \notin [\underline{b},\overline{b}] \enspace .
%\end{alignat*}
See \cite{Hofschuster2004} for details about how \texttt{C-XSC} manipulates rectangular containment sets over $\mathbb{IR}$ and $\mathbb{IC}$.

\section{Results} \label{SectionResults}

We have implemented the Sch\"onhage-Strassen algorithm, our containment-set version with rectangular complex intervals and long multiplication in \texttt{C++} using the \texttt{C-XSC} library \cite{Hofschuster2004}. Our implementation is available at \href{http://www.math.canterbury.ac.nz/~r.sainudiin/codes/capa/multiply/}{\url{http://www.math.canterbury.ac.nz/~r.sainudiin/codes/capa/multiply/}}
%\href{http://www.thomassteinke.org/home/multiply.cpp}{http://www.thomassteinke.org/home/multiply.cpp}
\footnote{Please also download the \texttt{C-XSC} library from \href{http://www.math.uni-wuppertal.de/~xsc/}{\url{http://www.math.uni-wuppertal.de/~xsc/}}.}. Results show that, using base $256$, our version of the algorithm is usually able to guarantee correct answers for up to 75,000-digit numbers.

The following graph compares the speed of long multiplication (labelled `Long multiplication'), the conventional Sch\"onhage-Strassen algorithm with different underlying data types (the implementation using the \texttt{C-XSC} \texttt{complex} data type is the line labelled `complex na\"ive SS' and the one using our own implementation of complex numbers based on the \texttt{C++} \texttt{double} data type is labelled `double na\"ive SS') and our containment-set version (`cinterval extended SS') on uniformly-random $n$-digit base-$256$ inputs. All tests were performed on a 2.2
GHz 64-bit AMD Athlon 3500+ Processor running Ubuntu 9.04 using \texttt{C-XSC} version 2.2.4 and gcc version 4.3.1. Times were recorded using the \texttt{C++} \texttt{clock()} function --- that is to say, CPU time was recorded.  Note that only the `Long multiplication' and `cinterval extended SS' implementations are guaranteed to produce correct results. The `double na\"ive SS' and `complex na\"ive SS' implementations may have produced erroneous results, as the implementations do not necessarily provide sufficient precision; these are still shown for comparison. Note also that by `na\"ive' we mean that these implementations use fixed-precision floating-point arithmetic, whereas the `real' Sch\"onhage-Strassen algorithm uses variable-precision, which is much slower.

\begin{center}
\epsfig{file=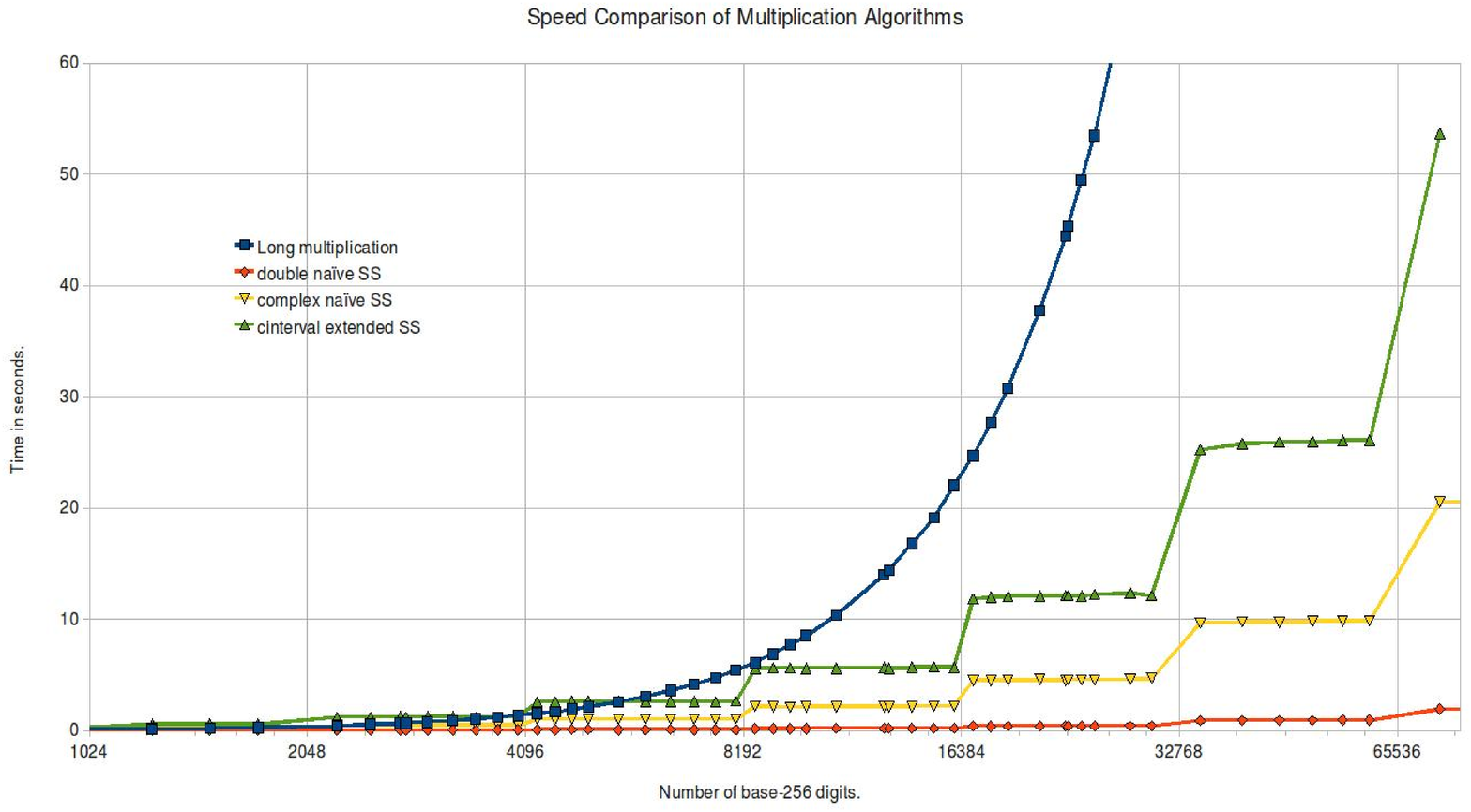, height=250pt}
\end{center}

The above graph shows that the Sch\"onhage-Strassen algorithm is much more efficient than long multiplication for large inputs. However, our modified version of the algorithm is slower than the na\"ive Sch\"onhage-Strassen algorithm. We believe that \texttt{C-XSC} is not well-optimised; for example, their punctual \texttt{complex} data type (used in the `complex na\"ive SS' implementation) is much slower than our \texttt{double}-based complex data type (used in the `double na\"ive SS' implementation), even though ostensibly they are the same thing. We see that the \texttt{C-XSC} \texttt{cinterval} type (used in the `cinterval extended SS' implementation) is about three times as slow as the \texttt{complex} type. This leaves the possibility that a more optimised implementation of containment sets would be able to compete with commercial algorithms. Investigations such as \cite{CFASTER} have shown that the Na\"ive Sch\"onhage-Strassen algorithm is able to compete with commercial implementations.

Note that the ``steps'' seen in the graph can be explained by the fact that the algorithm will always round the size up to the nearest power of two. Thus there are steps at the powers of two.  The most important feature of our results is the range of input sizes for which our algorithm successfully determines the answer. Using only standard double-precision IEEE floating-point numbers, we are able to use the algorithm to multiply 75,000-digit, or 600,000-bit, integers; this range is more than sufficient for most applications, and at this point the second Sch\"onhage-Strassen algorithm will become competitive.

\section{Conclusion}\label{Conclusion}

Our investigation has demonstrated that the Sch\"onhage-Strassen algorithm with containment sets is a practical algorithm that could be used reliably for applications such as cryptography. Sch\"ohage and Strassen never showed that their algorithm had any practical value. However, as our implementation was not optimised, this is more of a feasibility study than a finished product.

Note that the advantage of our algorithm over the original Sch\"onhage-Strassen algorithm is that we make use of hardware-based floating-point arithmetic, whereas the original is designed to make use of much slower software-based arithmetic. Both the original algorithm and our adaptation always produce correct results. However, we use a different approach to guaranteeing correctness. The na\"ive algorithms we mention are not guaranteed to be correct because they are modifications of the Sch\"onhage-Strassen algorithm which does not take measures to ensure correctness --- they simply use fixed-precision floating-point arithemtic and hope for the best; these are only useful for speed comparisons.

It remains to optimise our implementation of the algorithm to compete with commercial libraries. This is dependent on a faster implementation of interval arithmetic. It may also be interesting to use circular containment sets rather than rectangular containment sets. The advantage of circular containment sets is that they easily deal with complex rotations --- that is, multiplying by $e^{{\rm i} \theta}$; this is in fact the only type of complex multiplication (other than division by an integer) that our algorithm performs.

\bibliographystyle{eptcs} % or whatever you prefer

\begin{thebibliography}{1}
\providecommand{\bibitemstart}[1]{\bibitem{#1}}
\providecommand{\bibitemend}{}
\providecommand{\bibliographystart}{}
\providecommand{\bibliographyend}{}
\providecommand{\url}[1]{\texttt{#1}}
\providecommand{\urlprefix}{Available at }
\providecommand{\bibinfo}[2]{#2}
\bibliographystart

\bibitemstart{FUERER}
\bibinfo{author}{Martin F\"urer} (\bibinfo{year}{2007}):
  \emph{\bibinfo{title}{Faster Integer Multiplication}}.
\newblock In: {\sl \bibinfo{booktitle}{39th ACM STOC}}, \bibinfo{address}{San
  Diego, California, USA}, pp. \bibinfo{pages}{57--66}.
\bibitemend

\bibitemstart{CFASTER}
\bibinfo{author}{Pierrick Gaudry}, \bibinfo{author}{Alexander Kruppa} \&
  \bibinfo{author}{Paul Zimmermann} (\bibinfo{year}{2007}):
  \emph{\bibinfo{title}{A gmp-based implementation of sch\"{o}nhage-strassen's
  large integer multiplication algorithm}}.
\newblock In: {\sl \bibinfo{booktitle}{ISSAC '07: Proceedings of the 2007
  international symposium on Symbolic and algebraic computation}},
  \bibinfo{publisher}{ACM}, \bibinfo{address}{New York, NY, USA}, pp.
  \bibinfo{pages}{167--174}.
\bibitemend

\bibitemstart{Hofschuster2004}
\bibinfo{author}{Hofschuster} \& \bibinfo{author}{Kr\"amer}
  (\bibinfo{year}{2004}): \emph{\bibinfo{title}{{C-XSC 2.0: A {\tt C++}}
  library for extended scientific computing}}.
\newblock In: \bibinfo{editor}{R~Alt}, \bibinfo{editor}{A~Frommer},
  \bibinfo{editor}{RB~Kearfott} \& \bibinfo{editor}{W~Luther}, editors: {\sl
  \bibinfo{booktitle}{Numerical software with result verification}}, {\sl
  \bibinfo{series}{Lecture notes in computer science}} \bibinfo{volume}{2991},
  \bibinfo{publisher}{Springer-Verlag}, pp. \bibinfo{pages}{15--35}.
\bibitemend

\bibitemstart{KNUTH}
\bibinfo{author}{Donald~E. Knuth} (\bibinfo{year}{1998}):
  \emph{\bibinfo{title}{The Art of Computer Programming}},
  ~\bibinfo{volume}{2}.
\newblock \bibinfo{publisher}{Addison-Wesley}, \bibinfo{edition}{3} edition.
\bibitemend

\bibitemstart{Moore2009}
\bibinfo{author}{Ramon~E. Moore}, \bibinfo{author}{R.~Baker Kearfott} \&
  \bibinfo{author}{Michael~J. Cloud} (\bibinfo{year}{2009}):
  \emph{\bibinfo{title}{Introduction to Interval Analysis}}.
\newblock \bibinfo{publisher}{Society for Industrial and Applied Mathematics},
  \bibinfo{address}{Philadelphia, PA, USA}.
\bibitemend

\bibitemstart{ABOUND}
\bibinfo{author}{Colin Percival} (\bibinfo{year}{2003}):
  \emph{\bibinfo{title}{Rapid multiplication modulo the sum and difference of
  highly composite numbers}}.
\newblock {\sl \bibinfo{journal}{Math. Comput.}}
  \bibinfo{volume}{72}(\bibinfo{number}{241}), pp. \bibinfo{pages}{387--395}.
\bibitemend

\bibitemstart{RSA}
\bibinfo{author}{R.L. Rivest}, \bibinfo{author}{A.~Shamir} \&
  \bibinfo{author}{L.~Adleman} (\bibinfo{year}{1977}): \emph{\bibinfo{title}{A
  Method for Obtaining Digital Signatures and Public-Key Cryptosystems}}.
\newblock {\sl \bibinfo{journal}{Communications of the ACM}}
  \bibinfo{volume}{21}(\bibinfo{number}{2}), pp. \bibinfo{pages}{120--126}.
\bibitemend

\bibitemstart{SCHOENHAGE}
\bibinfo{author}{A.~Sch\"onhage} \& \bibinfo{author}{V.~Strassen}
  (\bibinfo{year}{1971}): \emph{\bibinfo{title}{Schnelle Multiplikation gro\ss
  er Zahlen (Fast Multiplication of Large Numbers)}}.
\newblock {\sl \bibinfo{journal}{Computing: Archiv f\"ur elektronisches Rechnen
  (Archives for electronic computing)}} \bibinfo{volume}{7}, pp.
  \bibinfo{pages}{281--292}.
\newblock \bibinfo{note}{(German)}.
\bibitemend

\bibliographyend
\end{thebibliography}
%\begin{thebibliography}{1}

%\bibitem{GA08}
%R.J.~van Glabbeek, C.~Author \& Y.S.~Else (2008):
%\newblock \emph{An example of a paper with a rather large title-to-content ratio}.
%\newblock {\sl Electronic Proceedings in Theoretical Computer Science} 0,
% pp. 1--3.
%\newblock Available at \url{http://style.eptcs.org/}.

%\end{thebibliography}

\end{document}